\documentclass[10pt,twosided]{article}

\usepackage[all]{xy}
\usepackage{amsmath}
\usepackage{amsthm}

\begin{document}
\input{defs.sty}

\newcommand{\pa}[1]{\frac{\partial}{\partial{#1}}}

\title{
On  parametrization of linear
pseudo-differential boundary value
control systems}

\author{Jouko  Tervo, Markku
Nihtil\"a and Petri Kokkonen
\thanks{Department of Mathematics and Statistics,
University of Kuopio, P.O.Box 1627,
FIN-70211 Kuopio, Finland. Email: {\tt Jouko.Tervo@uku.fi}}}

\maketitle

\begin{abstract}

The paper considers pseudo-differential boundary value control
systems.
The underlying operators form an algebra ${\s D}$
with the help of which we are able to formulate typical boundary value
control problems.
The symbolic calculus gives tools to form
e.g. compositions, formal adjoints, generalized right or left inverses
and compatibility conditions.
By a parametrizability we mean that for
a given control system ${\bf A}u=0$ one finds an operator
${\bf S}$ such that ${\bf A}u=0$ if and only if $u={\bf S}f$.
The computation rules of ${\s D}$  (or its appropriate subalgebra
${\s D}'$)
guarantee that in many applications ${\bf S}$ can be
refinely analyzed or even explicitly
calculated.
We outline some methods
of homological algebra
 for the study of parametrization ${\bf S}$.
Especially the projectivity of a certain
factor module (defined by the system equations) implies the
parametrizability.
We give some examples to illustrate our computational methods.

\end{abstract}
\vskip1cm

{\bf AMS-classification}. 93B25, 93C20, 93A10,  35S99

\section{Introduction}

The general theory of  control systems corresponding to
boundary value problems for linear and nonlinear partial
differential equations (PDEs) has been developed by two
apparently different approaches: One has applied
{\it functional analytic}
 and { \it algebraic} methods. 
In addition, {\it differential geometry} increasingly offers  effective
tools and gives geometric intuition in this field. 
Functional analytic methods are based
 on the analysis of underlying  operators defined
 abstractly in appropriate (Banach) spaces, e.g. 
 \cite{banks83,fattorini99,curtain95}. 
 Definitions of basic concepts such as
 controllability,
 observability, stability are 
 typically formulated using {\it trajectories}
 which have direct intuitions to the real world  problems
 in technology and in science.

 In frequence space problems (e.g. transfer function analysis)
 algebraic approach is more conventional 
 \cite{vidyasagar85,quadrat02} but in
 state space problems  algebraic study of control systems
 is quite recent method.
 The algebraic methods analyze {\it algebraic
 structures} such as  structures of certain modules,
 exact sequences, compatibility
 conditions and one-sided inverses.
In these approaches one mixes
the input, state and output variables
(which are called { control variables}) and studies only the
structural properties of this mixed system
(e.g. \cite{willems91,pommaret02,fliess95,fliess98}).
In the certain linear cases one
is able e.g. to apply the related
(differential) algebraic modules
to determine whether the system is controllable or not.
The algebraic methods have often
 less intuitive connections to real problems.

Here we consider {\it parametrizability} 
of certain linear control systems which means that for
a given control system ${\bf A}u=0$ one finds an operator
${\bf S}$ such that ${\bf A}u=0$ if and only if $u={\bf S}f$.
The operators ${\bf A}$ and ${\bf S}$
are built on the operators of ${\s D}$ where ${\s D}$ is the
chosen algebra of operators.
Parametrizability is a kind
of structural (internal) property
of the control system. 
The concept is related to the differential {\it flatness}
of (ordinary differential)
control systems
(\cite{fliess95,fliess98,chetverikov02,levine06}).  
The "parameter $f$" is in the role of
flat output. $f$ is {\it endogeneous} if it can be conversely
expressed with the help of $u$ by 
$
f=  {\bf Q}u
$
where ${\bf Q}$  is built on the operators of ${\s D}$.
In this case the system can be considered to be "flat".
In some cases parametrizability is  equivalent to
the {\it torsion freeness} of certain structural modules
(\cite{pommaret02}).

By now the structural study of state space problems, 
especially from the algebraic point of view, is
mainly developed without boundary values. Some attempts
nevertheless can
be found in literature \cite{fliess98,petit01,nihtila01}.
It is clear that the structural
properties related to the boundary value control problems
 depend also on the boundary conditions, not only the
PDE system.
In the following we suggest how to take into account the
boundary conditions in parametrization.
We use very general formulation of boundary value systems.
The operators consists of appropriate {\it pseudo-differential 
and boundary value operators}.
This calculus enables, for example, the consideration of 
the trace operator. In addition, 
the calculus gives the greater freedom
in the manipulation of e.g. the compositions, adjoints,
compatbility conditions and left/right inverses.
We express our formalism tailored for the boundary value problems
related to certain partial differential equations. 
Besides boundary value control systems for PDEs, our approach
enables also e.g. the consideration of delay systems.

Firstly
we choose an appropriate {\it algebra} ${\s D}$ of pseudo-differential
and boundary value operators.
This algebra
 originated in \cite{monvel71} and afterwards it has been
enlarged (\cite{grubb86,krainer00}).
The algebra gives very natural frame to calculate with
PDEs and related boundary value operators.
We generalize the definition of
intrinsic parametrization 
(the word "intrinsic" refers to the
situation where state, input and output
variables are mixed)
for the control systems
which are corresponding to 
 this algebra and we give some preliminary
tools to study and find the parametrization. 
Compared with
the earlier differential algebraic
methods we loose certain
 structures such as  differential
fields and the derivation rules typical in the context
of partial differential operators.
In our case the modules are not generally over entire rings
(integral domains).
This follows e.g. because of the ring $C^\infty(\ol G\times\Delta)$
is not entire: We have $v_1v_2=0$ for any nonzero functions whose
supports are disjoint.
Furthermore ${\s D}$ as the ring of matrices is not entire.

In section 4 we formulate some preliminary ideas to study
the parametrizability by certain homological algebraic methods.
We define (more general) ${\s D}'$-parametrizability concept, where
${\s D}'$ is a {\it subalgebra with unity of ${\s D}$}.
We don't try to make clear
here relations between (various variants of)
controllability and parametrizability.
The definition of structural controllability for the boundary value
control systems (under consideration) can be potentially
founded on the use of torsion freeness concept which can be
generalized also for nonintegral domains.
Because the underlying module is not generally over an entire
ring the application of torsion concept
is slightly  more complicated than in the conventional case.

One aim of this paper is to show that the applied algebra gives 
tools to construct the operators needed in parametrization.
Parametrizabilty gives potential 
methods for practical controller design
since from 
$u={\bf S}f$ it follows that we can express, say 
the input variables
$c=(c_1,...,c_p)$,
the state variables
$v=(v_1,...,v_q)$ and
the output variables
$y=(y_1,...,y_r)$ with the help of "(free) variables $f$" as
\be
c={\bf S}_1f,\ v={\bf S}_2f,\ y={\bf S}_3f.
\ee
Here the operators ${\bf S}_j$ are  
constitued of the operators lying in the algebra ${\s D}$
and so we can calculate $c,\ v,\ y$ using computational rules
of ${\s D}$.
For example, in the output tracking problem one seeks
the "parameter" $f=f^*$ such that (at least
optimally) $y^*={\bf S}_3f^*$
where $y^*$ is a given reference output. The required
input is simply $c=c^*={\bf S}_1f^*$.
This kind of "open controls" are important in some problems
of modern technology (in molecular physics, for example).
The other potential field of applications is the optimal
control problems. Roughly speaking this option
can be described as follows: If we have an optimal control problem
$\min_{u} F(u)$ under the contraint ${\bf A}u=0$, we can transform it
to unconstrained problem $\min_{f} F({\bf S}f)$.

\subsection{Basic notations}       \label{sec1.1}

Let $G$ be an open bounded set in $\R^n$ and let $\Delta$ be an interval
in $\R$.
We assume that the closure $\ol G$ is a
smooth differentiable manifold with boundary.
Furthermore
we suppose that the
boundary $\partial G$ is orientable and
that the unit normal vector $\nu=\nu(x),\ x\in \partial G$
is pointing "outwardly" on $\partial G$.

The spaces
$C^\infty(\ol G)$,
$C^\infty(\partial G)$,
$C^\infty(\ol G\times\Delta)$
and $C^\infty(\partial G\times\Delta)$
are correspondingly the collections
of smooth functions
$\ol G\to\R$,
$\partial G\to\R$,
$\ol G\times\Delta\to\R$
and
$\partial G\times\Delta\to\R$.
$C_c^\infty(G)$ is the space of test functions
and $S(\R^n)$ denotes the Schwartz space of rapidly decreasing
functions $\R^n\to\R$.
Furthermore, we define
\bean
C^\infty(G\times\ol\R_+)&=&\{\phi\in
C^\infty(G\times\R_+)\ |\ {\rm for\ which\ there}\\
&& {\rm exists\ a\ }
\psi\in C^\infty(G\times \R)\ {\rm such\ that\ }
\phi=\psi|_{G\times\ol\R_+}\}
\eean
and
\bean
C^\infty_c(G\times\ol\R_+)&=&\{\phi\in
C^\infty(G\times\ol\R_+)\ |\ \phi\ {\rm has\ a\ compact\ support}\}.
\eean
$D'(G)$ is the space of distributions that is, $D'(G)$ is the dual of
$C_c^\infty(G)$.

Denote
\ben
{\s R}^{N,m}=
C^\infty(\ol G\times\Delta)^N
\oplus C^\infty(\partial G\times\Delta)^m.
\een
In addition, we denote
\ben
{\s R}={\s R}^{1,1}=
C^\infty(\ol G\times\Delta)
\oplus C^\infty(\partial G\times\Delta).
\een
Let 
\bean
{\s R}_0^{N,m}
&=&\{(v,w)\in {\s R}^{N,m}|\ \textrm{$v$\ and $w$\ have\ compact\ supports}\
\\
&&
\textrm{in \ $\ol G\times\Delta$\ and in $\partial G\times\Delta$,\ respectively} \}.
\eean
For $(v,w)\in {\s R}^{N,m},\ (\phi,\psi)\in {\s R}_0^{N,m}$
we can define an inner product
\be
\la (v,w),(\phi,\psi)\ra=
\sum_{j=1}^N\la v_{j},\phi_j\ra_{L_2(G\times\Delta)}+
\sum_{k=1}^m\la w_k,\psi_k\ra_{L_2(\partial G\times\Delta)}.
\ee

Let $\hat u$ (or ${\s F}u$) denote the standard Fourier transform
\be
\hat u(\xi)=({\s F}u)(\xi)=\int_{\R^n}u(x)e^{-{\rm i}\la x,\xi\ra}
{\rm d}x
\ee
for a function $u$ in the Lebesgue space $L_1(\R^n)$.
The inverse Fourier transform is for appropriate 
functions
\ben
{\s F}^{-1}f(x)=(2\pi)^{-n}
\int_{\R^n}\hat f(\xi)e^{{\rm i}\la x,\xi\ra}{\rm d}\xi.
\een
The ring of $m\times n$-matrices
is denoted by $M(m\times n)$.
The partial derivative ${{\partial^\alpha}\over{\partial x^\alpha}}$
is denoted more shortly by $\partial_x^\alpha$.

\section{Control system}\label{sec2}

\subsection{Pseudo-differential and boundary
value operators}   \label{sec2.1}

We formulate the (generalized) boundary control systems 
 with the help of pseudo-differential 
 and boundary value operators
(e.g. \cite{grubb86,krainer00}).
Let $S^m_\nu(G \times\Delta\times\R^{n+1})$ be a space of functions,
so called anisotropic symbols,
$a\in C^\infty(G\times\Delta)$ such that
for any compact set $K\subset G\times\Delta,\ \alpha,\ \beta
\in \N_0^n,\ j,l\in \N_0$
there exists a constant      $C_{\alpha,\beta,j,l,K} >0$
such that
\be\label{2.1a}
|\partial_x^\alpha\partial_t^j
\partial_\xi^\beta\partial_\eta^l a(x,t,\xi,\eta)|
\leq C_{\alpha,\beta,j,l,K}
\rho(\xi,\eta)^{m-|\beta|-\nu l}\ {\rm for\ all}\ (\xi,\eta)
\in\R^{n+1},\ (x,t)\in K
\ee
where
\ben
\rho(\xi,\eta)=(1+|\xi|^2+|\eta|^{2/\nu})^{1/2}.
\een
As is standard, the space of symbols 
$S^m_\nu(G \times\Delta\times\R^{n+1})$ can be equipped with
the topology defined by the appropriate semi-norms.
In the case where we have no boundary considerations,
the (classical) pseudo-differential operator $A:C_c^\infty(G\times\Delta)
\to C^\infty(G\times\Delta)$ is defined by
\be\label{2.1b}
A\phi(x,t)=(2\pi)^{-(n+1)}\int_{\R^n}\int_{\R}
e^{{\rm i}\la (x,t),(\xi,\eta)\ra} a(x,t,\xi,\eta)\hat \phi
(\xi,\eta) {\rm d}\xi{\rm d}\eta.
\ee
When the boundary $\partial G$ is included in the considerations the
definition of pseudo-differential operator contains some
modifications.
Especially, the so called transmission condition is needed.
In addition, we need several classes of other operators
(boundary operators).

It is sufficient to explain the operator classes only in the 
local case that is, in the case
when the closure $\ol G$ is replaced with a set $G'\times\ol\R_+$ where
$G'$ is an open subset of $\R^{n-1}$ and
where the boundary $\partial
G$ is replaced with a set $G'$.
The technically tedious
reduction to that case is based on the partition of unity
and it is standard.

Generally speaking, the symbols of "boundary operators" satisfy conditions
analogous, in a sense, to the transmission condition. This
provides the inclusion
$$
\mathcal{A}
\bigl(C_c^{\infty}(G'\times\overline{\R}_+\times\Delta)
\oplus C_c^{\infty}(G'\times\Delta)\bigr)\subset
C^{\infty}(G'\times\overline{\R}_+\times\Delta) \oplus
C^{\infty}(G'\times\Delta)
$$
for any boundary value operator
$\mathcal{A}$ explained in section \ref{sec2.2.1}.
In addition, the kernels of operators involved in a boundary value problem
must have proper supports. Then $$ \mathcal{A}
\bigl(C^{\infty}(G'\times\overline{\R}_+\times\Delta)
\oplus C^{\infty}(G'\times\Delta)\bigr)\subset
C^{\infty}(G'\times\overline{\R}_+\times\Delta) \oplus
C^{\infty}(G'\times\Delta) $$ and one can compose the operators of
boundary problems.
In the following we explain the operators in detail.

{\bf A pseudodifferential operator}
($\psi$do)
$$
r^+A \colon
C_c^{\infty}(G'\times\overline{\R}_+\times\Delta)
\rightarrow D'(G'\times\R_+\times\Delta)
$$
is
defined by $$
r^+Av(x,t)=(2\pi)^{-(n+1)}r^+\int_{\R^n}\int_{\R}
e^{i\langle(x,t),(\xi,\eta)\rangle} a(x,t,\xi,\eta)
(\widehat{e^+v})(\xi,\eta)\,d\xi\,d\eta. $$
The symbol $a$ for any
compact $K\subset G'\times\overline{{\R}}_+\times\Delta$
satisfies
\begin{equation} \label{estim}
|\partial_x^{\alpha}\partial_t^j\partial_{\xi}^{\beta}\partial_{\eta}^l
a(x,t,\xi,\eta)|
\leq C_{\alpha,j,\beta,l,K}\,
\rho(\xi,\eta)^{m-|\beta|-\nu l},\quad (x,t,\xi,\eta)\in
K\times\R^{n+1},
\end{equation}
where 
$\alpha,\beta\in\N^n_0$, $j,l\in\N_0$. 
The space of symbols satisfying (\ref{estim}) is denoted by
$S^m_\nu(G'\times\ol\R_+\times\Delta\times\R^{n+1})$.
The number
$m$ is called the {\it order} of the symbol $a$ and, simultaneously, the
order of $\psi$do $r^+A$.
Above $r^+$ refers to the restriction operator 
$r^+f=f|_{ G'\times\R_+\times\Delta}$
and $e^+$ refers to the extension by zero from
$G'\times \ol\R_+\times\Delta$ on  $G'\times \R\times\Delta$ .

The symbol $a$ satisfies the {\it transmission condition}, if
the expansion
\begin{equation}              \label{trans}
\partial^{\gamma}_{x_n}a(x,t,\xi,\eta)\big|_{x_n=0}=
\sum_{p=0}^m\alpha_{\gamma p}(x',t,\xi',\eta)\xi_n^p+
\sum_{k=-\infty}^{\infty}a_{\gamma k}(x',t,\xi',\eta)
{\bigl(\rho(\xi',\eta)-i\xi_n\bigr)^k\over
\bigl(\rho(\xi',\eta)+i\xi_n\bigr)^{k+1}}
\end {equation}
holds for all $\gamma\in\N_0$ while $\alpha_{\gamma p}\in
S_{\nu}^{m-p}(G'\times\Delta\times\R^n)$, $a_{\gamma k}$ is rapidly
decreasing sequence in $S_{\nu}^{m+1}(G'\times\Delta\times\R^n)$
(that is, for any semi-norm $p$ on 
$S_{\nu}^{m+1}(G'\times\Delta\times\R^n)$ and any $N\in\N$ there
exists a constant $c_{N,p}$ such that $p(a_{\gamma k})\leq
c_{N,p}(1+k)^{-N}$)
and
$\rho(\xi',\eta)=(1+|\xi'|^2+|\eta|^{2/\nu})^{1/2}$. The space of
symbols of order $m$ satisfying the transmission condition will be
denoted by ${\bf U}_{\nu}^m$.

\noindent {\bf Assertion 1.} {\it If a $\psi$do $r^+A$ satisfies the transmission
condition then $r^+Av\in
C^{\infty}(G'\times\overline{\R}_+\times\Delta)$ for any
$v\in C_c^{\infty}(G'\times\overline{\R}_+\times\Delta)$.}

The consideration of Assertion 1
can be reduced to the one-dimensional case
in the same way as in \cite{monvel71}, Theorem 2.9.
Partial differential operators are typical examples of pseudo-differentials
operators.

{\bf A potential operator.} A function $k\in
C^{\infty}(G'\times\overline{\R}_+\times\Delta
\times\R^{n+1})$ is called a potential symbol of order $m$
if
\be\label{po}
 k(x,t,\xi',\xi_n,\eta)= \sum_{k=0}^{\infty} a_k
(x,t,\xi',\eta) \bigl(\rho(\xi',\eta)-i\xi_n\bigr)^k
\bigl(\rho(\xi',\eta)+i\xi_n\bigr)^{-(k+1)},
\ee
where $a_k$ is a
rapidly decreasing sequence in the space
$S_{\nu}^{m+1}(G'\times\overline{\R}_+\times\Delta
\times\R^n)$ consisting of functions subjected to
(\ref{estim}) with $\rho(\xi,\eta)$ replaced by $\rho(\xi',\eta)$. The
space of potential symbols of order $m$ is denoted by ${\bf
K}^m_{\nu}$.

A potential operator is defined by
$$ Kw(x,t)=
(2\pi)^{(n+1)}\int_{\R^n}\int_{\R}
e^{i\langle(x,t),(\xi,\eta)\rangle} k(x,t,\xi',\xi_n,\eta)
\hat{w}(\xi',\eta)\,d\xi\,d\eta $$
for $w\in
C_c^{\infty}(G'\times\Delta)$. It is obvious that
$K\bigl(C_c^{\infty}(G'\times\Delta)\bigr)\subset
C^{\infty}(G'\times\overline{\R}_+\times\Delta)$.
The (classical) solution operators of certain partial differential
boundary value problems contain integral terms
\be
Kg(x,t)=\int_{\partial G\times\Delta} 
G(x,t,y',\tau)g(y',\tau) {\rm d}\sigma {\rm d}\tau
\ee
which are examples of  potential operators. The other
examples are the adjoints
of the below defined trace operators which are often potential operators.

{\bf A trace operator.} A function $t\in
C^{\infty}(G'\times\Delta\times\R^{n+1})$ is called a
trace symbol of {\it order} $m$ and {\it class} $d$ if
\begin{equation}          \label{trace}
t(x',t,\xi',\xi_n,\eta)=
\sum_{p=0}^{d-1}\alpha_p(x',t,\xi',\eta)\xi_n^p+
\sum_{k=0}^{\infty} a_k (x',t,\xi',\eta)
{\bigl(\rho(\xi',\eta)+i\xi_n\bigr)^k\over
\bigl(\rho(\xi',\eta)-i\xi_n\bigr)^{k+1}}
\end{equation}
where $\alpha_p\in S_{\nu}^{m-p}(G'\times\Delta\times\R^n)$
and $a_k$ is a rapidly decreasing sequence in the space
$S_{\nu}^{m+1}(G'\times\Delta\times\R^n)$.
The space of trace symbols of order $m$ and class $d$ is denoted by
${\bf T}^{m,d}_{\nu}$.

A trace operator is defined by
$$ Tv(x',t)=
(2\pi)^{-n}\int_{\R^{n-1}}\int_{\R}
e^{i\langle(x',t),(\xi',\eta)\rangle} \int^+ t(x',t,\xi',\xi_n,\eta)
(\widehat{e^
+v})(\xi,\eta)\,d\xi_n\,d\xi'\,d\eta $$
for $v\in
C_c^{\infty}(G'\times\overline{\R}_+\times\Delta)$. 
Here $\int^+fd\xi_n$ is the curve integral $\int_\gamma f(z){\rm d}z$
where $\gamma$ is "a path rounding the upper complex half plane
$\mathrm{Im} z>0$ in the counterclockwise direction" that is,
$\int^+fd\xi_n= \lim_{R\to\infty}\oint_{\Gamma_R}f(z)dz$ where
$\Gamma_R$ is the boundary of the half ball $B(0,R)\cap
\{z\in \C|\ \mathrm{Im} z>0\}$.
The
operator $T$ can be written in the form
$Tv=\sum_{p=0}^{d-1}S_p\bigl(\partial^p_{x_n}v(x',0,t)\bigr)+T_0v$
where $S_p$ are $\psi$dos
with symbols in $S^{m-p}_{\nu}(G'\times\Delta\times\R^n)$
and $T_0$ is of class $0$.
It is obvious that
$T\bigl(C_c^{\infty}(G'\times\overline{\R}_+\times\Delta)\bigr)
\subset C^{\infty}(G'\times\Delta)$.
Restrictions of partial differential operators on the boundary
(which usually appear in boundary conditions) are typical
trace operators.

{\bf A Green operator.} A function $b\in C^{\infty}
(G'\times\overline{\R}_+\times\Delta\times\R^{n+2})$
is called a Green symbol of {\it order} $m$ and {\it class} $d$ if
\begin{eqnarray}\label{go}
&&b(x,t,\xi',\xi_n,\zeta_n,\eta)
\\
&&=\sum_{p=0}^{d-1}k_p(x,t,\xi',\xi_n,\eta)\zeta_n^p+
\sum_{j,l=0}^{\infty} a_{jl}(x',t,\xi',\eta)
{(\rho(\xi',\eta)-i\xi_n)^j\over (\rho(\xi',\eta)+i\xi_n)^{j+1}}
{(\rho(\xi',\eta)+i\zeta_n)^l\over (\rho(\xi',\eta)-i\zeta_n)^{l+1}},\nonumber
\end{eqnarray}
where $k_p\in {\bf K}_{\nu}^{m-p}$ and $a_{jl}$ is a rapidly
decreasing double sequence in the space
$S_{\nu}^{m+2}(G'\times\Delta\times\R^n)$. The space of Green symbols of
order $m$ and class $d$ is denoted by ${\bf B}^{m,d}_{\nu}$.

A Green operator is defined by
$$ Bv(x,t)=
(2\pi)^{-n-1}\int_{\R^{n}}\int_{\R}
e^{i\langle(x,t),(\xi,\eta)\rangle} \int^+ b(x,t,\xi',\xi_n,\zeta_n,\eta)
(\widehat{e^+v})(\xi',\zeta_n,\eta)\,d\zeta_n\,d\xi\,d\eta
$$ 
for
$v\in C_c^{\infty}(G'\times\overline{\R}_+\times\Delta)$.
The operator $B$ can be written in the form
$Bv=\sum_{p=0}^{d-1}K_p\bigl(\partial^p_{x_n}v(x',0,t)\bigr)+B_0v$
where $K_p$
are potential operators with
symbols in ${\bf K}^{m-p}_{\nu}$ and $B_0$ is of class
$0$. It is obvious that
$B\bigl(C_c^{\infty}(G'\times\overline{\R}_+\times\Delta)\bigr)\subset
C^{\infty}(G'\times\overline{\R}_+\times\Delta)$.
Typical sources of singular Green operators are the compositions of
potential and trace operators. They are also born  in the compositions
of truncated pseudo-differential operators. Finally they appear 
e.g. in many
solution operators related to the usual boundary value problems.

{\bf A pseudodifferential operator on the boundary}
$$
Q \colon
C_c^{\infty}(G'\times\Delta)
\rightarrow C^\infty(G'\times\Delta)
$$
is
defined by 
$$
Qw(x',t)=(2\pi)^{-n}\int_{\R^{n-1}}\int_{\R}
e^{i\langle(x',t),(\xi',\eta)\rangle} q(x',t,\xi',\eta)
\hat w(\xi',\eta) {\rm d}\xi'{\rm d}\eta. $$ 
with
the symbol $q\in S^m_\nu(G'\times \Delta\times\R^n)$.
Partial differential operators on the boundary (manifold)
form an example of these operators.

In practice the operator $r^+A$ is often
a partial differential operator.
The operator $T$
is a natural generalization of the usual partial differential trace
operator appeared in classical theory of boundary value
problems. The operators $K$, $B$ and $Q$ are  needed e.g. in the
consideration of the compositions, inverses and adjoints.
The operators $T,\ K,\ B,\  Q$  are called
{\it boundary operators}.
Together with these operators 
one can formulate very rich variety of boundary value
problems. This more general formulation also  gives 
more symmetry in the calculus of adjoints and  compositions.

\begin{example}    \label{ex1}

Let $L(D)=a\partial_x^2+b\partial_x+c$ be a PDO with
constant coefficients and let $d_1,\ d_2$ be constants
(in this example  we have no time variable).
Furthermore, let $G=]0,1[\subset\R$. Then
under relevant assumptions the solution $v_1$ for   the system
\bea
&&L(D)v_1=q(x)v_2                       \\
&&\partial_xv_1(0)+d_1v_1(0)=0,\ \partial_xv_1(1)+d_2v_1(1)=0
\eea
is formally given by  (\cite{nihtila01})
\be
v=
{ P}_0 (qv_2)+{ R}_0(qv_2)
\ee
where   $P_0$ is a pseudo-differential operator
\be
P_0 \phi (x)=(2\pi)^{-1}\int_{\R}
e^{{\rm i}\la x,\xi\ra}{1\over{L(\xi)}}\hat \phi (\xi)
{\rm d}\xi
\ee
and  $R_0$ is  a singular Green operator
(note that in one-dimensional case the integration reduces to
summation)
\bea
R_0 \phi (x)&=&
-{{g_1(x,0)}\over {W(0)}}[({1\over 2}+g_2(0,0))P_0\phi (1)
-g_2(0,1)P_0 \phi (0)]\nonumber\\
&&-
{{g_2(x,0)}\over {W(0)}}[g_1(0,0))P_0\phi (1)
-({1\over 2}+g_1(0,1))P_0 \phi (0)].
\eea
Here $L(\xi)=-a\xi^2+b{\rm i}\xi-c$ 
is the {\it symbol} of $L(D)$ and 
\ben
g_k(x,\lambda)=(2\pi)^{-1}{\bf pv}\int_{\R}
\left[e^{{\rm i}\xi x}{{g_k(\xi)}\over{\lambda-L(\xi)}}
{\rm d}\xi\right],\ \lambda=0,1,\ k=1,2,
\een
\ben
W(0)=g_1(0,0)g_2(1,0)-(1+g_1(1,0))(1+g_2(0,0)),
\een
where
\ben
g_1(\xi)=(ad_2-b-a{\rm i}\xi)e^{-{\rm i}\xi},
\
g_2(\xi)=b-ad_1+a{\rm i}\xi .
\een
Here ${\bf pv}$  denotes that the integral is taken
in the sense of  principal value.
\end{example}

\subsection{Computation rules}\label{sec2.2}

We survey some basic 
computational properties of the above defined operators.

\subsubsection{Operator algebra}            \label{sec2.2.1}

{\bf A boundary value operator}
$\mathcal{A}=\qmatrix{
r^+A+B & K \cr T & Q\cr}$
is called {\it proper} if the kernels of all operators $r^+A, B, K, T, Q$
have proper supports.
The space of proper boundary value operators
is denoted by $\mathcal{D}$

\noindent {\bf Assertion 2.} {\it If $\mathcal{A}$ and $\mathcal{B}$ are in
$\mathcal{D}$, then $\mathcal{A}\mathcal{B}
:={\s A}\circ {\s B}\in\mathcal{D}$.}

To verify this assertion one can modify the arguments used
in \cite{grubb86} 
for the case of parameter dependent boundary value problems.

We say that the operator $\mathcal{A}=\qmatrix{
r^+A+B & K \cr T & Q \cr}$ is of class
$\mathcal{D}^{m,d}_{\nu}$,
$$
m=\qmatrix{ m_1,m_2 & m_3
\cr m_4 & m_5 \cr}, \qquad d=\qmatrix{ d_2 \cr
d_4 \cr}
$$
when

the symbol of  $r^+A$ is in the class  ${\bf U}^{m_1}_{\nu}$,

the symbol of $B$ is in the class ${\bf B}^{m_2,d_2}_{\nu}$,

the symbol of $K$ is in the class  ${\bf K}^{m_3}_{\nu}$,

the symbol of $T$ is in the class ${\bf T}^{m_4,d_4}_{\nu}$,

the symbol of $Q$ is in the class $S^{m_5}_{\nu}$.

The assertion 2 implies that
proper boundary value operators of arbitrary order $m$ and
class $d$ form an {\it algebra} $\mathcal{D}$ with respect to
standard addition and composition of operators.
In addition,
from the above definitions and assertions it follows
that the spaces $\mathcal{D}$ and $
C^{\infty}(G'\times\overline{\R}_+\times\Delta)\oplus
C^{\infty}(G'\times\Delta)$ are $\mathcal{D}$-{\it modules}.
The order and the class of ${\s AB}$ can be calculated
with the help of orders and classes of ${\s A}$ and ${\s B}$.

In the global setting we have operators
\ben
r^+A:C^\infty(\ol G\times\Delta)\to C^\infty(\ol G\times\Delta)
\een
\ben
K:C^\infty(\partial G\times\Delta)\to C^\infty(\ol G\times\Delta)
\een
\ben
T:C^\infty(\ol G\times\Delta)\to C^\infty(\partial G\times\Delta)
\een
\ben
B:C^\infty(\ol G\times\Delta)\to C^\infty(\ol G\times\Delta)
\een
\ben
Q:C^\infty(\partial G\times\Delta)\to C^\infty(\partial G\times\Delta)
\een
and
\ben
{\s A}:{\s R}\to {\s R}
\een
where ${\s R}=
C^\infty(\ol G\times\Delta)\oplus C^\infty(\partial G\times\Delta)$.
In addition ${\s R}$ is a ${\s D}$-module.

\subsubsection{Formal adjoint}      \label{sec2.2.2}

As mentioned above in 
${\s R}=C^\infty(\ol G\times\Delta) 
\oplus C^\infty(\partial G\times\Delta)$ we use the inner product 
(e.g. for functions $U_1\in {\s R},\ U_2\in {\s R}_0$)
\ben
\la U_1,U_2\ra=\la v_{1},  v_{2}\ra_{L_2(G\times\Delta)}+
\la w_{1}, w_{2}\ra_{L_2(\partial G\times\Delta)}
\een
where
\ben
U_k= (v_{k},w_{k}),\ k=1,2.
\een

For any boundary value operator
$\mathcal{A}\in\mathcal{D}$ there exists an operator $$
\mathcal{B}:
C_c^{\infty}(\ol G\times\Delta)\oplus
C_c^{\infty}(\partial G\times\Delta)\rightarrow
{D}'(G\times\Delta)\oplus
{D}'(\partial G\times\Delta) $$ such that
$(\mathcal{A}U,V)=( U,\mathcal{B}V)$
where $(\cdot,\cdot)$ is the standard duality form on
\ben
(C_c^{\infty}(\ol G\times\Delta)\oplus
C_c^{\infty}(\partial G\times\Delta))\times
({D}'(G\times\Delta)\oplus
{D}'(\partial G\times\Delta)).
\een
If the operator $\mathcal{B}$
maps the space $\mathcal{R}$ into itself we will say that formal
adjoint to $\mathcal{A}$ exists and will denote the operator
$\mathcal{B}$ by $\mathcal{A}^*$.
When the operator ${\s A}^*:
{\s R}\to {\s R}$ 
exists we have
\be                         \label{4.11}
\la V,{\s A}U\ra=   \la{\s A}^*V,U\ra
\ee
for $V,U\in {\s R}_0$.

\noindent {\bf Assertion 3}. {\it The adjoint
operator $\mathcal{A}^*$ exists and is in the algebra
$\mathcal{D}$ if and only if the operator $\mathcal{A}=
\qmatrix{
r^+A+B & K \cr T & Q\cr}$ is of class $0$ (i.e., the
class of  operators $B$ and $T$ equals $0$) and in the
transmission condition (\ref{trans}) for the symbol of $r^+A$ all
$\alpha_{\gamma p}$, except $\alpha_{\gamma 0}$, are equal to $0$.
In the case where
$\mathcal{A}^{*}$ exists then
 there exists $\mathcal{A}^{**}$ and $\mathcal{A}^{**}=\mathcal{A}$.}

We omit here the proof of Assertion 3 but we remark that analogous
results and considerations can be found in  the
monograph \cite{grubb86}.

We find that for
an operator $\mathcal{A}\in\mathcal{D}$ the formal adjoint
$\mathcal{A}^*$, generally speaking, does not exist.
However, one can consider subset $\mathcal{D}_1\subset\mathcal{D}$ that
consists of boundary value operators satisfying conditions of
the Assertion 3. Then $\mathcal{D}_1$ is an algebra and
for any $\mathcal{A}\in\mathcal{D}_1$ there exists the
formal adjoint $\mathcal{A}^*$ in $\mathcal{D}_1$.
If we fix $m_0=\qmatrix{ 0,-1 & -1/2 \cr -1/2 & 0
\cr }$ and $d_0=\qmatrix{ 0 \cr 0 \cr}$, then the subalgebra
$\mathcal{D}^{m_0,d_0}_{\nu}\subset\mathcal{D}_1$ possesses the
same properties as $\mathcal{D}_1$.

\begin{example}               \label{ex2}

Let $G=]0,1[\subset\R$ and
\ben
Tv=d(\partial)v|_{\partial G}
\een
where $d(\partial)$ is a pseudo-differential operator
\ben
d(\partial)v (x)=(2\pi)^{-1}\int_{\R}
e^{{\rm i}\xi x}d(\xi)\widehat{e^+v}(\xi)
{\rm d}\xi,
\een
such that its symbol
$d(\xi)$ satisfies
\be\label{te}
|d(\xi)|\leq C{1\over{(1+|\xi|)^\kappa}},\ \kappa>1.
\ee
Then we have  for $v\in C^\infty(\ol G),\  g\in C^\infty(\partial G)$
\bean
&&\la Tv,g\ra_{L_2(\partial G)}
=
(Tv)(0)\ol{g(0)}+ (Tv)(1)\ol{g(1)}\\
&&=
(2\pi)^{-1}\ol{g(0)}\int_{\R}d(\xi)\widehat{e^+ v}(\xi){\rm d}\xi
+
(2\pi)^{-1}
\ol{ g(1)}\int_{\R}e^{{\rm i}\xi}d(\xi)\widehat{e^+ v}(\xi){\rm d}\xi\\
&&=
(2\pi)^{-1}\left
(\int_G\left(\ol{g(0)}\int_{\R}e^{-{\rm i}\xi x}d(\xi){\rm d}\xi\right)
v(x){\rm d}x
+
\int_G\left(\ol{g(1)}\int_{\R}e^{{\rm i}\xi(-x+1)} d(\xi)
{\rm d}\xi\right)v(x){\rm d}x\right)\\
&&=
\la v, T^*g\ra_{L_2(G)}
\eean
where
\be\label{tad}
T^*g(x)=
(2\pi)^{-1}\left(g(0)\int_{\R}e^{-{\rm i}\xi x}\ol{d(\xi)}{\rm d}\xi
+
g(1)\int_{\R}e^{{\rm i}\xi(-x+1)}\ol{d(\xi)}
{\rm d}\xi\right).
\ee
We see that $T^*$ is a potential-type operator, but it does not
generally possess the transmission property. Consideration
of operators like (\ref{tad}) can be found in \cite{vishik67}.
\end{example}

\subsubsection{Order reducing operators}         \label{sec2.2.3}

The so called {\it order
reduction} can be used to help the forming of adjoints
in parametrization processes.

Let $\mathcal{A}$ be a boundary value operator of order
$m=\qmatrix{ m_1,m_2 & m_3 \cr m_4 &
m_5\cr}$ and class $d=\qmatrix{ d_2 \cr
d_4 \cr}$. If $T=\qmatrix{
r^+T_1 & 0 \cr 0 & T_2 \cr}$ while
$\mbox{ord}\,r^+T_1=-N_1$
(that is, the order of $r^+T_1$ is $-N_1$)
and $\mbox{ord}\,T_2=-N_2$,
$N_1,\ N_2\in\N$ then the
operator $\mathcal{A}\circ T$ is of order
$\qmatrix{ m_1-N_1,m_2-N_1 & m_3-N_2 \cr m_4-N_1 &
m_5-N_2 \cr}$ and class $\qmatrix{
d_2-N_1 \cr d_4-N_1 \cr}$.
Note that the adjoint operator $(\mathcal{A}\circ T)^*$
exists if $N_1>\max\,(m_1,d_2,d_4)$.

\subsection{System equations}

In the following
we denote more shortly by ${\s A}=\qmatrix{{\s A}_1&{\s A}_2\cr
{\s A}_3&{\s A}_4\cr}$
a typical element 
${\s A}=\qmatrix{r^+A+B& K\cr T&Q\cr}$ of
${\s D}$.

Let $v_1,...,v_N$ and $w_1,...,w_m$ be indeterminates in
$C^\infty(\ol G\times\Delta)$ and in
$C^\infty(\partial G\times\Delta)$, respectively.
Furthermore, let
$({\s A}_{1ij}),\ ({\s A}_{2ik}),\ ({\s A}_{3lj}),
\ ({\s A}_{4lk})$ 
be operator matrices where
${\s A}_{1ij}$ are of type
$r^+A_{ij}+B_{ij}$,
${\s A}_{2ik}$ are of type
$K_{ik}$,
${\s A}_{3lj}$ are of type
$T_{lj}$ and
${\s A}_{4lk}$ are of type
$Q_{lk}$.
We consider
 the control system
\bea        \label{2.12a}
\sum_{j=1}^N{\s A}_{1ij}v_j+
\sum_{k=1}^m{\s A}_{2ik}w_k&=&0,\quad i=1,...,N_1,\\
 \sum_{j=1}^N {\s A}_{3lj}v_j+
\sum_{k=1}^m{\s A}_{4lk}w_k&=&0,\quad l=1,...,m_1.\nonumber
\eea
The state, input and output variables are not separated.

The system (\ref{2.12a}) can be given in a matrix form as follows.
Let  
$({\s A}_{1ij})\in M(N_1\times N)$,
$ ({\s A}_{2ik})\in M(N_1\times m)$,\
$({\s A}_{3lj})\in M(m_1\times N)$ and
$({\s A}_{4lk})\in M(m_1\times m)$ be the matrices corresponding to
the system (\ref{2.12a}). 
Denote
$
u=\qmatrix{v\cr w\cr}$
where $v=\qmatrix{v_1\cr\vdots\cr v_N\cr},\ 
w=\qmatrix{w_1\cr\vdots\cr w_m\cr}$ and denote
\be        \label{of}
 {\bf A}=\qmatrix{({\s A}_{1ij})&({\s A}_{2ik})\cr
 ({\s A}_{3lj})&
({\s A}_{4lk})\cr}\in M((N_1+m_1)\times(N+m)).
\ee
Then the system (\ref{2.12a}) is equivalent to the equation
\be           \label{2.13}
 {\bf A}u=0.
\ee
The operator ${\bf A}$ is a linear operator $
{\s R}^{N,m}
\to {\s R}^{N_1,m_1}.$  We say the the operator $\bf A$ is in
${\s D}(N_1+m_1,N+m)$ when it is of the form
(\ref{of}).
This notation uniquelly indicates the types of submatrices.
Using the standard
{\it Frechet space} topologies in $C^\infty(\ol G\times\Delta)$
and in
$C^\infty(\partial G\times\Delta)$ one has

\noindent {\bf Assertion 4}. 
{\it The linear operator ${\bf A}:{\s R}^{N,m}
\to {\s R}^{N_1,m_1}$ is continuous.}

\begin{example}    \label{ex3}

In this example we outline how to
choose the above operators in the case
of second order linear PDEs.
Consider a linear 
 coupled system of partial differential equations
\be\label{2.14}
\sum_{j=1}^N A_{ij}(x,t,\partial)v_j
=0,\ 
\ee
for $i=1,...,N_1$, where $A_{ij}(x,t,\partial)$
are, for example, the second order operators
\ben
A_{ij}(x,t,\partial)=
\sum_{|\sigma|\leq 2}a_\sigma^{ij}(x,t)\partial^\sigma,
\een
where
\ben
\partial=\left({\partial\over{\partial
x_1}},...,
{\partial\over{\partial x_n}}, {\partial\over{\partial t}}\right)
=:(\partial_x,\partial_t).
\een
Denote ${\s A}_{1ij}(x,t,\partial)=A_{ij}(x,t,\partial)$.

We assume that the solution satisfies the following
 homogeneous boundary
 conditions
\be\label{2.16}
\left[\sum_{j=1}^N d_{lj}(x,t,\partial_x)v_j
\right]\vert_{\partial G\times\Delta}=0
 ,\ l=1,...,m_1
\ee
where $d_{lj}(x,t,\partial_x)$ are
first order partial differential operators
\ben
d_{lj}(x,t,\partial_x)=\sum_{k=1}^nd_{ljk}(x,t)
{{\partial}\over{\partial {x_k}}}
+d_{lj0}(x,t).
\een

In this example we choose the operator matrix $({\s A}_{1ij})$ such that
\ben           \label{2.17}
({\s A}_{1ij})=
\qmatrix{r^+A_{11}(x,t,\partial)&\cdots&
r^+A_{1N}(x,t,\partial))\cr\vdots\cr
r^+A_{N_1 1}(x,t,\partial))&\cdots&
r^+A_{N_1 N}(x,t,\partial))\cr},
\een
where $r^+f:=f|_{G\times\Delta}$ is (as above) the
restriction operator on $G\times\Delta$. The boundary operators $B_{ij},\ K_{ik},
\ Q_{lk}$ are
zero operators and the operators
${\s A}_{3lj}$ corresponding to the boundary operators
$T_{lj}$ are defined by
\be             \label{2.18}
({\s A}_{3lj})
=
\qmatrix{r'd_{11}(x,t,\partial_x)&\cdots&
r'd_{1N}(x,t,\partial_x)\cr\vdots\cr
r'd_{m_1 1}(x,t,\partial_x)&\cdots&
r'd_{ m_1 N}(x,t,\partial_x)\cr}
\ee
 where
$r'f:=f|_{\partial G\times\Delta}$  is the {\rm restriction operator}
on $\partial G\times\Delta$.
As a conclusion we find that
\be\label{2.19}
 {\bf A}=\qmatrix{({\s A}_{1ij})&0\cr ({\s A}_{3lj})&0\cr}.
\ee
\end{example}

\begin{remark}
When we choose $d_{lj}(x,t,\partial_x)=0,\ 
 l=1,...,m_1,\ j=1,...,N$ the function $v_j$ does not
satisfy any boundary condition on $\partial G$.

\end{remark}

\begin{example}\label{ex4}
Consider the delay system
\be         \label{delay}
{{dv}\over{dt}}=A_0v+\sum_{k=1}^p A_kv(t-h_k)+B_0V
\ee
where $0<h_1<\cdots < h_p$, $v(t)\in \C^N,\  V(t)\in\C^q$ and
$A_0,\ A_k\in M(N\times N),\ B_0\in M(N\times q)$.

Noting that
\bea\label{dela}
&&f(t-h_k)=(2\pi)^{-1}\int_{\R}\hat f(\eta)e^{{\rm i}(t-h_k)\eta}{ \rm d}
\eta \\
&&=
(2\pi)^{-1}\int_{\R}e^{-{\rm i}h_k\eta}\hat f(\eta)e^{{\rm i}t\eta}{ \rm d}
\eta
=
P_kf(t)\nonumber
\eea
we see that    the system (\ref{delay}) can be written in the form
\be
{\bf A}\qmatrix{v\cr V\cr}=0
\ee
where
\be
{\bf A}= \qmatrix{{d\over{dt}}-A_0-\sum_{k=1}^pA_kP_k & -B_0\cr}.
\ee
Here $P_k$ is a pseudodiffererential operator with symbol
$p_k(\eta)=e^{-{\rm i}h_k\eta}$.
\end{example}

\section{Parametrization of a control system applying
algebraic methods}\label{sec3}

\subsection{Parametrizability of a control system}\label{sec3.1}

Here we consider parametrizability of the
above defined boundary value control
systems.
The control system
\bea        \label{4.4}
\sum_{j=1}^N{\s A}_{1ij}v_j+\sum_{k=1}^m{\s A}_{2ik}w_k&=&0,\ i=1,...,N_1,\\
 \sum_{j=1}^N {\s A}_{3lj}v_j+
\sum_{j=1}^m{\s A}_{4lk}w_k&=&0,\ l=1,...,m_1,\nonumber
\eea
or more simply
\be\label{4.5}
{\bf A}u=0 
\ee
is said to be {\it parametrizable}, if there exists
$N',\ m'\in\N$ and a linear operator
\be  \label{4.7a}
{\bf S}=
\qmatrix{({\s S}_{1jp})&
({\s S}_{2jq})
\cr ({\s S}_{3kp}) &
({\s S}_{4kq})\cr} \in {\s D}(N+m, N'+m')
\ee
such that 
\be\label{4.5a}
{\bf A}u=0\Leftrightarrow u={\bf S}f
\ee
(the components $f_m$ of $f$ are ${\s D}$-linear independent).
The condition (\ref{4.5a}) means that the equation
${\bf A}u=0$ determinates all {\it compatibility  conditions}
of the equation ${\bf S}f=u$.

\begin{example}
Let $G\subset\R^2$. Denote $\partial_1=\pa{x_1}$, $\partial_2=\pa{x_2}$
and let $I$, $I'$ be the identity
mappings on $C^\infty(\ol{G})$ and $C^\infty(\partial G)$, respectively.
Consider a system
\ben
\partial_2 v_1-\partial_1 v_1-\partial_2 v_2+v_2=0.
\een
The system can be put into the
form (the inclusion of $w_1$ is due to the notational convenience)
\ben
\qmatrix{\partial_2-\partial_1 & -
\partial_2+I & 0 \cr 0 & 0 & 0}\qmatrix{v_1 \cr v_2 \cr w_1}=0.
\een
The system has a parametrization (\cite{pommaret02})
\be
\qmatrix{v_1\cr v_2\cr w_1\cr}=
\qmatrix{\partial_2-I&0\cr \partial_2-\partial_1&0\cr 0&I'\cr}
\qmatrix{f\cr w_1\cr}=:{\bf S} \qmatrix{f\cr w_1\cr}.
\ee
\end{example}

\begin{example}    \label{ex1a}
Consider the control system related to Example \ref{ex1}
\bea
&&L(D)v_1=q(x)v_2                      \\
&&\partial_xv_1(0)+d_1v_1(0)=0,\ \partial_xv_1(1)+d_2v_1(1)=0\nonumber\\
&&y=v_1|_{\partial G}.
\eea
This system can be put into the form
\be\label{ex5:eq2}
\qmatrix{L(D)&-q(x)I&0\cr T&0&0\cr r'&0&-I'\cr}
\qmatrix{
v_1\cr v_2\cr y\cr}=0.
\ee
Here $I$ (resp. $I'$) is the identity mapping on $C^\infty(\ol G)$
(resp. $C^\infty(\partial G)$). $T$ is the trace operator
\be
Tv_1=\begin{cases}
\partial_xv_1(x)+d_1v_1(x),\ & x=0\cr
\partial_xv_1(x)+d_2v_1(x),\ & x=1\cr
\end{cases}.
\ee
Due to Example \ref{ex1} the
 system (\ref{ex5:eq2}) has a parametrization given by
\be
\qmatrix{
v_1\cr v_2\cr y\cr}=\qmatrix{(P_0+R_0)(q(\cdot))\cr I\cr
r'(P_0+R_0)(q(\cdot))\cr}v_2
\ee
\end{example}

The homological algebra gives  tools to study
structurally the
parametrizability. We describe the basic idea as follows.
Suppose that ${\s D}$ is an algebra of suitable operators.
Let $P_1,\ P_2,\ M$ be ${\s D}$-modules
and let $d_1:P_2\to P_1$ and $d_2:P_1\to M$
be ${\s D}$-homomorphisms.
Recall that the sequence of modules
\ben
\begin{xy}
\xymatrix{
\ar[r] P_2\ar[r]^{d_1} & P_1\ar[r]^{d_2} & M\ar[r] & 0
}
\end{xy}
\een
is   {\it complex} if ${\rm im}\ d_1\subset {\rm ker}\ d_2$
 and it is {\it exact} if
${\rm ker}\ d_2={\rm im}\ d_1$.
The exactness means that $d_2v=0$ if and only if $v=d_1f$
and so the solutions of system $d_2v=0$ can be {\it
parametrized} by $d_1$.
The {\it homology} of the above complex is  the quotient
${\rm ker}\ d_2/
{\rm im}\ d_1$.  The homology  measures how much the sequence
differs from being exact. Furthermore, 
using certain homology groups, 
Ext$^n(M,{A})$, one can eventually
specify how far a given module is e.g. from being  projective. 
Using the groups Tor$^n(M,{A})$ one is able to
measure how far a given module is from being flat.

\subsection{System modules}\label{sec4.1}

As we found above, ${\s D}$ is an algebra of operators. Hence
it is a (left) ${\s D}$-module
(for basic concepts of homological algebra
see e.g. \cite{rotman79,kunz85,osborne99,northcott76}).
The module ${\s D}$ has a unit 
\be
{\s I}=\qmatrix{r^+I&0\cr 0&I'\cr}
\ee
where $I$ is the identity operator
$C^\infty(\ol G\times\Delta)\to 
 C^\infty(\ol G\times\Delta)$ and 
 $I'$ is the identity operator
$C^\infty(\partial G\times\Delta)\to
 C^\infty(\partial G\times\Delta)$.
 Hence ${\s D}$ is a unitary ${\s D}$-module.

We at first {\it choose a subalgebra }
 ${\s D}'$ of ${\s D}$ which contains an unit $I$.
 Furthermore, let
  ${\s D}'^N$ be the direct product of modules
\be
{\s D}'^N={\s D}'\times \cdots \times {\s D}'.
\ee
Then ${\s D}'^N$ is a (left) free  ${\s D}'$-module whose canonical basis is
\be
{E}_1=({\s I},0,\cdots,0),\ ...\ {E}_N=(0,\cdots,0,{\s I}).
\ee

\begin{example} \label{ex7}
The following sets of operators are subalgebras of ${\s D}$
\begin{enumerate}
\item $\displaystyle {\s D}'=\{\qmatrix{r^+A(D)&0\cr 0&0\cr}|
\ A(D)\ {\rm is\ a\ PDO\ with\ constant\
coefficients}\}$,

\item $\displaystyle{\s D}'=\{\qmatrix{r^+A(x,t,D)
&0\cr 0&0\cr}|\ A(x,t,D)\ {\rm is\ a\ PDO\ with\ real\ analytic\
coefficients}\}$,

\item $\displaystyle {\s D}'=\{\qmatrix{r^+A(x,t,D)&0\cr r'd(x,t,D)
& aI'\cr}|\ A(x,t,D),\ d(x,t,D)\ \textrm{
are\ PDOs\ with\ }
C^\infty(\ol G\times
\Delta)\textrm{-\ coefficients},\ a\ \in\R\}$,

\item $\displaystyle {\s D}'={\s D}_1$,
where ${\s D}_1$ is as in subsection \ref{sec2.2.2}.
\end{enumerate}

\end{example}

Suppose that the control system
\bea        \label{4.1}
\sum_{j=1}^N{\s A}_{1ij}v_j+
\sum_{k=1}^m{\s A}_{2ik}w_k&=&0,\ i=1,...,N_1,\\
 \sum_{j=1}^N {\s A}_{3lj}v_j+
\sum_{k=1}^m{\s A}_{4lk}w_k&=&0,\ l=1,...,m_1,\nonumber
\eea
is given.
Let $\ol N=\max\{N,m\},\ \ol N_1=\max\{N_1,m_1\}$.
The system (\ref{4.1}) can be expressed equivalently in the form
\be                                                    \label{sl}
\qmatrix{{\s L}
_{11}&\cdots&{\s L}_{1\ol N}\cr \cdots&\cdots&\cdots\cr
{\s L}_{\ol N_11}&\cdots&{\s L}_{\ol N_1\ol N}\cr}
\qmatrix{\qmatrix{v_1\cr w_1\cr }\cr \vdots\cr\qmatrix
{v_{\ol N}\cr w_{\ol N}\cr }\cr}=0.
\ee
For example,
suppose that $N>m,\ N_1> m_1$. Define
\be                             \label{3.2}
{\s L}_{ij}=\begin{cases}
\qmatrix{{\s A}_{1ij}&{\s A}_{2ij}\cr {\s A}_{3ij}&{\s A}_{4ij}\cr}
& \textrm{for}\ 1\leq j\leq m,\ 1\leq i\leq m_1,\cr
\qmatrix{{\s A}_{1ij}&0\cr {\s A}_{3ij}&0\cr}
& \textrm{for}\ m+1\leq j\leq N,\ 1\leq i\leq m_1,\cr
\qmatrix{{\s A}_{1ij}&{\s A}_{2ij}\cr 0&0\cr}
& \textrm{for}\ 1\leq j\leq m,\ m_1+1\leq i\leq N_1,\cr
\qmatrix{{\s A}_{1ij}&0\cr 0&0\cr}
& \textrm{for}\ m+1\leq j\leq N,\ m_1+1\leq i\leq N_1.\cr
\cr
\end{cases}
\ee
Similarly we can reformulate the other cases.

\begin{example}
The system
\ben
{\bf A}u=\qmatrix{{\s A}_{111}&{\s A}_{112} &{\s A}_{113}&
{\s A}_{211}&{\s A}_{212}\cr
{\s A}_{121}&{\s A}_{122} &{\s A}_{123}&{\s A}_{221}&{\s A}_{222}\cr
{\s A}_{311}&{\s A}_{312} &{\s A}_{313}&{\s A}_{411}&{\s A}_{412}\cr }
\qmatrix{v_1\cr v_2\cr v_3\cr w_1\cr w_2\cr}
=0
\een
can be put into the form
\ben
{\s L}
\qmatrix{\qmatrix{v_1\cr w_1\cr}\cr\qmatrix{ 
v_2\cr w_2\cr}\cr \qmatrix{ v_3\cr w_3\cr}\cr}=0 ,
\een
where
\ben
{\s L}=
\qmatrix{{\s L}_{11}&{\s L}_{12}& 
{\s L}_{13}\cr
{\s L}_{21}&{\s L}_{22}&{\s L}_{23}\cr}
=
\qmatrix{\qmatrix{{\s A}_{111}&{\s A}_{211}\cr
{\s A}_{311}&{\s A}_{411}\cr}&
\qmatrix{{\s A}_{112}&{\s A}_{212}\cr 
{\s A}_{312}&{\s A}_{412}\cr}&
\qmatrix{{\s A}_{113}&0\cr 
{\s A}_{313}&0\cr}\cr
\qmatrix{{\s A}_{121}&{\s A}_{221}\cr 
0&0\cr}&
\qmatrix{{\s A}_{122}&{\s A}_{222}\cr
0&0\cr}&
\qmatrix{{\s A}_{123}&0\cr 
0&0\cr}\cr}.
\een

\end{example}

Let ${\s D}'(m\times n)$ be the set of $m\times n$-matrices whose
elements are in ${\s D}'$.
Denote ${\s L}=({\s L}_{ij})\in {\s D}'(\ol N_1\times\ol N)$ and denote
$u_j=\qmatrix{v_j\cr w_j\cr},\ j=1,...,\ol N$. Let
$u=\qmatrix{u_1\cr\vdots\cr u_{\ol N}\cr}$. Then the system
(\ref{sl}) is
\be
{\s L}u=0.
\ee
Define a mapping $.{\s L}:{\s D}'^{\ol N_1}\to {\s D}'^{\ol N}$ by
\be
.{\s L}D=D{\s L}=\qmatrix{D_1&\cdots& D_{\ol N_1}\cr}
\qmatrix{{\s L}_{11}&\cdots&{\s L}_{1\ol N}\cr \cdots&\cdots&\cdots\cr
{\s L}_{\ol N_11}&\cdots&{\s L}_{\ol N_1\ol N}\cr}
\ee
for $D=  \qmatrix{D_1&\cdots& D_{\ol N_1}\cr}\in {\s D}'^{\ol N_1}$.
Then we find that $.{\s L}$ is a (left) ${\s D}'$-homomorphism that is,
\ben
.{\s L}(L_1D+L_2D')=L_1.{\s L}D+L_2.{\s L}D'
\een
for $L_1,\ L_2\in {\s D}',\ D,\ D'\in {\s D}'^{\ol N_1}$.
The image im$(.{\s L})$ is exactly ${\s D}'^{\ol N_1}{\s L}$ and it is
a submodule of ${\s D}'^{\ol N}$.
Hence we are able to define the (left) factor module
\be
M={\rm coker}(.{\s L})={\s D}'^{\ol N}/{\s D}'^{\ol N_1}{\s L}.
\ee
The elements of $M$ are of the form $[D]=D+ {\s D}'^{\ol N_1}{\s L}$.
Let $\pi:{\s D}'^{\ol N}\to M$ be the canonical surjection
\be
\pi D=[D].
\ee
$\pi$ is a surjective ${\s D}'$-homomorphism
and ${\rm ker}\ \pi ={\s D}'^{\ol N_1}{\s L}$ and so we have an
exact sequence
\be    \label{ra}
\xymatrix{
\ar[r]{\s D}'^{\ol N_1}\ar[r]^{.{\s L}}
& {\s D}'^{\ol N}\ar[r]^{\pi} & M\ar[r] & 0.  
}
\ee
The module $M$ corresponds to the control system (\ref{4.1}) in 
a sense that
\be
\sum_{j=1}^{\ol N}{\s L}_{ij}[E_j]=0,\ 1\leq i\leq \ol N_1,
\ee
that is, the system equations are always valid in $M$ for $[E_j]$.

\begin{remark}
In the case where the algebra ${\s D}'$ is commutative we are also able to
define a module $M':={\s D}'^{\ol N_1}/{\s L}{\s D}'^{\ol N}$.
In this case we find that
\be
{\s L}[D]=0,\ \forall D\in M'
\ee
that is, the system equations are valid in $M'$. Furthermore, we
have an exact sequence
\be       \label{ex}
\begin{xy}
\xymatrix{
\cdots \ar[r] & {\s D}'^{\ol N} \ar[r]^{{\s L}} & {\s D}'^{\ol N_1}\ar[r]& M' \ar[r] & 0
}.
\end{xy}
\ee
A simple example of the commutative subalgebra is the first case in
Example \ref{ex7}. The commutativity is very restrictive property
for ${\s D}'$.
The module $M'$ can be also defined in the case where ${\s D}'$ is only
an {\em Ore algebra}. In this case (\ref{ex}) is valid. We do not
consider module $M'$ here.

\end{remark}

The module $M$ is finitely generated because for every element\\
$[D]=[\qmatrix{D_1&\cdots&D_{\ol N}\cr}]\in M$
\be 
[D]=
\qmatrix{D_1&\cdots&D_{\ol N}\cr}+{\s D}'^{\ol N_1}{\s L}
=
\sum_{j=1}^{\ol N}D_j[E_j].
\ee
The elements $[E_j]$ do not, however, form
necessarily a basis
that is, every element $[D]\in M$
cannot be {\it uniquelly} expressed in the form
$[D]=\sum_{j=1}^{\ol N}D_j[E_j]$.
By (\ref{ra}) the module $M$ is {\it finitely presented}.

A ${\s D}'$-module $M$ is {\it free} if it has a basis $[F_1],...,[F_r]$.
Free module is always ${\s D}'$-isomorphic to the module ${\s D}'^r$.
A module $M$ is {\it projective
free} if there exists a ${\s D}'$-module $P$ such that the direct sum
$M\oplus P=F$  is free. 
The {\it  exact sequence} of projective ${\s D}'$-modules $P_j$
\be           \label{4.2}
\cdots  \xymatrix{\ar[r] & P_{n}\ar[r]^{d_n} & P_{n-1}\ar[r]&}\cdots \\
\xymatrix{\ar[r] & P_1\ar[r]^{d_1} & P_0\ar[r]^{d_0} & M\ar[r] & 0}
\ee
is a {\it projective resolution} of $M$. The projective resolution
always exists that is, we can choose   projective
${\s D}'$-modules $P_0,\ P_1,\ P_2,...$ and ${\s D}'$-homomorphisms
$d_0,\ d_1,\ d_2,...$ such that
\bean
d_0:P_0\to M\ {\rm is\ onto}\\
d_1:P_1\to {\rm ker}\ d_0\ {\rm is\ onto}\\
d_2:P_2\to {\rm ker}\ d_1\ {\rm is\ onto}\\
\eean
and so on.

\subsection{Algebraic criterions for parametrizability}\label{sec4.2}

Let
${\bf S}:R^{N',m'}\to R^{N,m}$ be an operator as in section \ref{sec3.1}
\be  \label{1}
{\bf S}=
\qmatrix{({\s S}_{1jp}')&
({\s S}_{2jq}')
\cr ({\s S}_{3kp}') &
({\s S}_{4kq}')\cr} \in {\s D}(N+m,N'+m')
\ee
Similarly to the construction
of ${\s L}$ define the matrix ${\s S}$ corresponding to ${\bf S}$ 
\be\label{sec4.2:eq2}
{\s S}=
\qmatrix{{\s S}_{11}&\cdots&{\s S}_{1\ol N'}\cr \cdots&\cdots&\cdots\cr
{\s S}_{\ol N1}&\cdots&{\s S}_{\ol N\ol N'}\cr}\in
{\s D}'(\ol N\times\ol N').
\ee
The homomorphism $.{\s S}:{\s D}'^{\ol N}\to {\s D}'^{\ol N'}$ is also
defined as $.{\s L}$.

Let $A$ be a left ${\s D}'$-
module. In this algebraic setting
we say,
more generally,
that the system ${\s L}X=0$,
$X\in A^{\ol{N}}$, is {\it ${\s D}'$-parametrizable
in $A^{\ol N}$}
if there exists a matrix ${\s S}$ given by (\ref{sec4.2:eq2}) such that
\be\label{ap}
{\s L}X=0\Leftrightarrow X=(X_1,...,X_{\ol N})=
{\s S}X'\ {\rm for}\ X'\in A^{\ol N'}.
\ee
In our case $A$ may be ${\s R}$ or ${\s D}'$, for example.
When $A={\s R}$ we get the parametrization concept in the sense
of subsection \ref{sec3.1}.

Denote Hom$(P,A)=\{d|\ d:P\to A\ {\rm is}\ {\s D}'-{\rm 
homomorphism}\}$
for ${\s D}'$-modules $P$ and $A$.
Hom$(P,A)$ is a group with respect to addition but
since ${\s D}'$ is not
necessarily commutative Hom$(P,A)$ is not generally a ${\s D}'$-module.
We find that the mapping
$
{\rm Hom}(.{\s L},{A}):{\rm Hom}({\s D}'^{\ol N},{A})
\to
{\rm Hom}({\s D}'^{\ol N_1},{A})$ defined by
\be
{\rm Hom}(.{\s L},{A})\phi=\phi\circ .{\s L}
\ee
is a ${\s D}'$-homomorphism.

\begin{lemma}       \label{le1}
The homomorphism ${\rm Hom}(.{\s L},{A})$ can be calculated by
\be
[{\rm Hom}(.{\s L},{A})\phi](D) =
\sum_{i=1}^{\ol N_1}
 D_i\left(\sum_{j=1}^{\ol N}{\s L}_{ij}\phi(E_j)\right)
\ee
for $D=(D_1,...,D_{\ol N_1})\in {\s D}'^{\ol N_1}$.
\end{lemma}

\begin{proof}
For any $\ol D=\sum_{j=1}^{\ol N}\ol D_j E_j\in
{\s D}'^{\ol N}$
and for  any $\phi\in {\rm Hom}
({\s D}'^{\ol N},{A})$ one has the expression
\ben
\phi(\ol D)=\sum_{j=1}^{\ol N}\ol D_j\phi(E_j)
\een
and so
\bea\label{ph}
[{\rm Hom}(.{\s L},{A})\phi](D) &=& (\phi\circ .{\s L})( D)=
\phi( .{\s L} D)= \phi( D{\s L})\\
&=&
\sum_{j=1}^{\ol N}( D{\s L})_j\phi(E_j)
=
\sum_{i=1}^{\ol N_1}
 D_i\left(\sum_{j=1}^{\ol N}{\s L}_{ij}\phi(E_j)\right)
\eea
for $D\in {\s D}'^{\ol N_1}$ as desired.
\end{proof}

\indent\begin{theorem}                          \label{th3}
The parametrizability condition
\be\label{apl}
{\s L}X=0\Leftrightarrow X={\s S}X'\ {\rm for}\ X'\in {A}^{\ol N'}.
\ee
is equivalent to the equality
\be        \label{apk}
{\rm ker\ Hom}(.{\s L},{A})={\rm im\ Hom}(.{\s S},{A}).
\ee
\end{theorem}

\begin{proof}
A. Suppose that (\ref{apl}) holds.
Then by Lemma \ref{le1} and by (\ref{apl}) we find that
\bea\label{p1}
&&
\phi\in {\rm ker\ Hom}(.{\s L},{A})\Leftrightarrow
{\rm  Hom}(.{\s L},{A})\phi=0\Leftrightarrow
\sum_{j=1}^{\ol N}{\s L}_{ij}\phi(E_j)=0
\nonumber\\
&&
\Leftrightarrow
\qmatrix{\phi(E_1)\cr\vdots\cr \phi(E_{\ol N})\cr}={\s S}X'
\Leftrightarrow
\phi(E_j)=\sum_{k=1}^{\ol N'}{\s S}_{jk}X_k'
\nonumber
\\
&&
\Leftrightarrow
\phi(D)=\sum_{j=1}^{\ol N}D_j\phi (E_j)
=
\sum_{j=1}^{\ol N}D_j\sum_{k=1}^{\ol N'}{\s S}_{jk}X_k'.
\eea
Define $\psi_{X'}\in {\rm Hom}({\s D}'^{\ol N'},A)$ by
$\psi_{X'}(\tilde D):=\sum_{k=1}^{\ol N'}\tilde D_kX_k'$.
Then we have $X_k'=\psi_{X'}(E_k)$ and so by Lemma \ref{le1}
\be
\phi(D)
=\sum_{j=1}^{\ol N}D_j\sum_{k=1}^{\ol N'}{\s S}_{jk}\psi_{X'}(E_k)
=
[{\rm Hom}(.{\s S},A)\psi_{X'}](D)
\Leftrightarrow \phi=        {\rm Hom}(.{\s S},A)\psi_{X'}.
\ee
Since for any $\psi\in {\rm Hom}({\s D}^{\ol N'},A)$
\be\label{es}
\psi=\psi_{X'},\ {\rm for}\ X'=\qmatrix{\psi(E_1)\cr\vdots\cr
\psi(E_{\ol N'})\cr}
\ee
the first part of the assertion is proved.

B. Conversely, suppose that (\ref{apk}) holds.
Define $\phi_X(D)=\sum_{j=1}^{\ol N} D_jX_j$. Then by Lemma 1
\bea\label{p2}
&&{\s L}X=0\Leftrightarrow \sum_{j=1}^{\ol N}{\s L}_{ij}X_j=0,\ \forall i
\Leftrightarrow
\sum_{j=1}^{\ol N}{\s L}_{ij}\phi_X(E_j)=0,\ \forall i\\
&&\Leftrightarrow
\sum_{i=1}^{\ol N_1}D_i\sum_{j=1}^{\ol N}{\s L}_{ij}\phi_X(E_j)=0,\ \forall D\in {\s D}'^{\ol{N}_1}
\Leftrightarrow
{\rm Hom}(.{\s L},A)\phi_X=0.\nonumber
\eea
By (\ref{apk}) the equivalence (\ref{p2}) is valid if and only if
\bea
&&
\phi_X={\rm Hom}(.{\s S},A)\psi
\Leftrightarrow
\phi_X( D)=[{\rm Hom}(.{\s S},A)\psi]( D)
=\sum_{j=1}^{\ol N} D_j
\sum_{k=1}^{\ol N'}{\s S}_{jk}\psi(E_k)
,\ \forall D\in {\s D}'^{\ol N}
\nonumber\\
&&
\Leftrightarrow
X_j=
\sum_{k=1}^{\ol N'}{\s S}_{jk}\psi(E_k)
\Leftrightarrow
X={\s S}X',\ {\rm for}\ X'=\qmatrix{\psi(E_1)\cr\vdots\cr \psi(E_{\ol N'})
\cr}.
\eea
Due to (\ref{es})
this completes the proof.
\end{proof}

The homomorphism $.{\s S}:{\s D}'^{\ol N}\to {\s D}'^{\ol N'}$
is similarly
defined as $.{\s L}$. Let $\tilde N$ be the ${\s D}'$-module
\be
\tilde N={\s D}'^{\ol N'}/{\s D}'^{\ol N}{\s S}.
\ee
For the definition and for the
basic properties of the homology
groups
${\rm Ext}^n(M,A)$
and ${\rm Tor}^n(M,A)$ (see e.g. \cite{rotman79,kunz85,osborne99})

We have

\indent\begin{theorem} \label{th:4.2.1}
Suppose that there exists a matrix (\ref{sec4.2:eq2}) such that
\be\label{eq:4.2.1}
D{\s L}=D'\Leftrightarrow D'{\s S}=0.
\ee
Then
the parametrizability condition (\ref{apk})
is valid if and only if
\be           \label{exb}
{\rm Ext}^1(\tilde N,A)=0.
\ee

\end{theorem}

\begin{proof}
One sees that the assumption (\ref{eq:4.2.1})
is equivalent to
\be                        \label{eq:4.2.2}
.{\s L}D=D'\Leftrightarrow .{\s S}D'=0.
\ee
Hence assumption  (\ref{eq:4.2.1}) is equivalent to the exactness
of the following sequence
\be\label{eq:4.2.3}
\begin{xy}
\xymatrix{
\dots\ar[r] &
{\s D}'^{\ol N_1} \ar[r]^{.{\s L}} &
{\s D}'^{\ol N} \ar[r]^{.{\s S}} &
{\s D}'^{\ol N'}
\ar[r]^{\tilde\pi } & \tilde N\ar[r] & 0    .
}
\end{xy}
\ee
The resolution (\ref{eq:4.2.3}) implies the truncated dual
sequence
\be\label{res3}
\begin{xy}
\xymatrix{
0\ar[r] & {\rm Hom}({\s D}'^{\ol N'},A)\ar[r]^{{\rm Hom}(.{\s S},A)} &
{\rm Hom}({\s D}'^{\ol N},A)\ar[r]^{{\rm Hom}(.{\s L},A)} &
{\rm Hom}({\s D}'^{\ol N_1},A) \ar[r] & \dots
}.
\end{xy}
\ee
Ext$^1(\tilde N,A):=
{\rm ker\ Hom}(.{\s L},{A})/{\rm im\ Hom}(.{\s S},{A})=0$
 if and only if the
 condition
(\ref{apk}) is valid. Hence
Theorem \ref{th3} implies the assertion.
\end{proof}

The condition (\ref{exb}) can also be
studied without any explicit calculations. For example, since
the module $\tilde N$ is finitely presented, the condition
(\ref{exb}) is valid when the module $\tilde N$ is flat.
In the case where $A$ is injective  Ext$^1(\tilde N,A)=0$.
Especially,
${\rm Ext}^1(\tilde N,{A})=0$ in the case where $\tilde N$ is projective.
Note that neither ${\s D}'$ nor ${\s R}$ are  generally
injective modules.

The conditions like
(\ref{exb}) are interesting because they are independent
of the used (projective) free resolutions. 
So they give valuable information for the instrinsic structure of
the system.

\begin{example}      \label{ex:4.2.1}

Let $G=]0,1[$ and $\Delta$ is any
open interval of $\R$ . Consider the system
\bean
{\q {v_1}t}&=&{\q {v_1}x}+v_2\\
 v_1(0,t)&=& v_1(1,t)=0.
\eean
In practice $v_1$ may be the state variable and $v_2$ may be the control
variable.
We have $N=2, \ N_1=1,\ m_1=1$ and 
$\partial G=\{0,1\}$. 
For notational convenience we may assume that $m=1$ although
no boundary value functions $w_1$ exist.
Hence $\ol N=2,\ \ol N_1=1$.
Denote $\partial_1={\partial\over{\partial x}},\
\partial_2={\partial\over{\partial t}}$.

We see that
\be
{\s L}:=
\qmatrix{{\s L}_1&{\s L}_2\cr}=
\qmatrix{\qmatrix{r^+(\partial_2^2-\partial_1^2)&0\cr r'I&0\cr}
&\qmatrix{-I&0\cr 0&0\cr}\cr}.
\ee

We seek a possible operator ${\s S}$ for which the assumption
(\ref{eq:4.2.1})
holds.
Let $D=\qmatrix{{\s A}_1&{\s A}_2\cr {\s A}_3&{\s A}_4\cr}\in {\s D}$ and
$(D',D'')\in {\s D}^2$
where $D'=\qmatrix{{\s A}_1'&{\s A}_2'\cr {\s A}_3'&{\s A}_4'\cr},\
D''=\qmatrix{{\s A}_1''&{\s A}_2''\cr {\s A}_3''&{\s A}_4''\cr}$.
We find that
\be\label{eq:4.2.4}
D{\s L}=D'
\ee
if and only if
\be           \label{eq:4.2.4a}
\qmatrix{\qmatrix{{\s A}_1r^+(\partial_2^2-\partial_1^2)+{\s A}_2r'I
&0\cr {\s A}_3r^+(\partial_2^2-\partial_1^2)+{\s A}_4r'I&0\cr}
&\qmatrix{-{\s A}_1&0\cr -{\s A}_3&0\cr}\cr}
=
\qmatrix{\qmatrix{{\s A}_1'&{\s A}_2'\cr {\s A}_3'&{\s A}_4'\cr} &
\qmatrix{{\s A}_1''&{\s A}_2''\cr {\s A}_3''&{\s A}_4''\cr}\cr}
\ee
if and only if
\bea\label{eq:4.2.4b}
&&
{\s A}_2''={\s A}_4''={\s A}_2'={\s A}_4'=0,\ -{\s A}_1={\s A}_1'',\
-{\s A}_3={\s A}_3'' \nonumber\\
&&
-{\s A}_1''r^+(\partial_2^2-\partial_1^2)+
{\s A}_2r'I={\s A}_1',\
-{\s A}_3''r^+(\partial_2^2-\partial_1^2)+
{\s A}_4r'I={\s A}_3.
\eea

Let $K$ be the potential operator
\be
Kg(x,t)=  xg(1,t)+(1-x)g(0,t).
\ee
Then we  observe that $r'Kg=g$ or $r'K=I'$.
Multiplying the equation
\be\label{eq:4.2.5}
-{\s A}_1''r^+(\partial_2^2-\partial_1^2)+
{\s A}_2r'I={\s A}_1'
\ee
by $K$ we get
\be\label{eq:4.2.6}
{\s A}_2={\s A}_1'K+{\s A}_1''r^+(\partial_2^2-\partial_1^2)K.
\ee
Similarly we get from the equation
\be\label{eq:4.2.7}
-{\s A}_3''r^+(\partial_2^2-\partial_1^2)+
{\s A}_4r'I={\s A}_3'
\ee
that
\be\label{eq:4.2.8}
{\s A}_4={\s A}_3'K+{\s A}_3''r^+(\partial_2^2-\partial_1^2)K.
\ee
Hence by (\ref{eq:4.2.4b})
\be\label{eq:4.2.9}
-{\s A}_1''r^+(\partial_2^2-\partial_1^2)+
({\s A}_1'K+{\s A}_1''r^+(\partial_2^2-\partial_1^2)K)r'I={\s A}_1'
\ee
and
\be\label{eq:4.2.10}
-{\s A}_3''r^+(\partial_2^2-\partial_1^2)+
({\s A}_3'K+{\s A}_3''r^+(\partial_2^2-\partial_1^2)K)r'I={\s A}_3'.
\ee
Expressing the equations (\ref{eq:4.2.9}), (\ref{eq:4.2.10})
and
\be\label{4.2.11}
{\s A}_2''={\s A}_4''={\s A}_2'={\s A}_4'=0
\ee
in   the matrix form we have
\be  \label{eq:4.2.12}
\qmatrix{\qmatrix{{\s A}_1'&{\s A}_2'\cr {\s A}_3'&{\s A}_4'\cr} &
\qmatrix{{\s A}_1''&{\s A}_2''\cr {\s A}_3''&{\s A}_4''\cr}\cr}{\s S}=0
\ee
where
\be\label{eq:4.2.13}
{\s S} =
\qmatrix{ \qmatrix{I_{1,0}-Kr'I
&0\cr 0&I\cr}&\qmatrix{0&0\cr 0&0\cr}\cr
\qmatrix{(r^+(\partial_2^2-\partial_1^2)
(-I+Kr'I)&0\cr 0&0\cr}&
\qmatrix{0&0\cr 0&I'\cr}\cr}.
\ee
In  section \ref{sec4.4} we shall find that ${\s S}$ is the
${\s D}'$-parametrization in any module $A^2$.
\end{example}

\subsection{Application of adjoints}\label{sec4.3}

The adjoints can  be applied
also in the algebraic analysis which we describe
in the sequel.
In this subsection we {\it assume that the required formal adjoints exist
in the subalgebra ${\s D}'$}.

\begin{remark}         \label{re:3.1.1}
If necessary we can apply order reducing to quarantee the existence
of formal adjoints. For example, in the case where $A={\s R}$ this can
be described as follows.

Suppose that ${\s T}$ and ${\s T}'$ are
isomorphisms ${\s R}^{\ol N}\rightarrow {\s R}^{\ol N}$ and
${\s R}^{\ol N_1}\rightarrow {\s R}^{\ol N_1}$, respectively. Then the
control system ${\s L}u=0$ is equivalent to the system
\be
({\s T}'\circ {\s L}\circ {\s T})u'=0
\ee
where ${\s T}u'=u$. Especially, the operators ${\s T}$ and ${\s T}'$
can be chosen to be order reducing that is, the orders and/or classes
of the composed operators are less than the original operators
(section \ref{sec2.2.3}). In parametrization we are able to use order
reducing operators to quarantee the existence of adjoints.
This is based on the following observation
\be
{\s L}u=0 \Leftrightarrow ({\s T}'\circ {\s L}\circ {\s T})u'=0
\Leftrightarrow u'={\s S}'f\Leftrightarrow u=({\s T}\circ
{\s S}')f\Leftrightarrow u={\s S}f
\ee
where ${\s S}'$ is the parametrization of the system
$({\s T}'\circ {\s L}\circ {\s T})u'=0$ and ${\s S}=
{\s T}\circ {\s S}'$.
\end{remark}

Let ${\s L}^*$ be the formal adjoint of ${\s L}$
\be
{\s L}^*=
\qmatrix{{\s L}_{11}^*&\cdots&{\s L}_{\ol N_1 1}^*\cr \cdots&\cdots&\cdots\cr
{\s L}_{1\ol N}^*&\cdots&{\s L}_{\ol N_1\ol N}^*\cr}
\in {\s D}'(\ol N_1\times \ol N)
\ee
defined as in \ref{sec2.2.2}.
The homomorphism $.{\s L}^*:{\s D}'^{\ol N}\to {\s D}'^{\ol N_1}$ 
is    similarly
defined as $.{\s L}$. Let $N$ be the ${\s D}'$-module
\be
N={\s D}'^{\ol N_1}/{\s D}'^{\ol N}{\s L}^*.
\ee
Suppose that $N$ has a {\it free resolution with finitely generated modules}
\be\label{res2}
\xymatrix{
\dots \ar[r] & {\s D}'^{\ol N'}\ar[r]^{.{\s S}^*} &
{\s D}'^{\ol N}\ar[r]^{.{\s L}^*} & {\s D}'^{\ol N_1}\ar[r]^{\pi^*} & N\ar[r] & 0.}
\ee
A ${\s D}'$-module $E$ is {\it stably free} if one can find
a finitely generated free module $F$ such that $E\oplus F$ is a finitely 
generated free
module that is, $E\oplus F$ is isomorphic to ${\s D}'^{n}$ for some $n$.
A projective stably free module $N$ always admits the free resolution
(\ref{res2}) with finitely generated modules, for example (\cite{lang93}).
In addition, in the case where ${\s D}'$ is Noetherian $N$ has the
free resolution with finitely generated  modules (\ref{res2}).
For example, the subalgebra
\ben
{\s D}'=\{\qmatrix{r^+A(D)&0\cr 0&0\cr}|
\ A(D)\ {\rm is\ a\ PDO\ with\ constant\
coefficients}\}.
\een
 is a commutative Noetherian integral domain.
Besides this kind of simple cases, substantial further work must be done
to analyze algebraic properties of more general systems.

\begin{remark}
{\rm Suppose that $A={\s D}'$}.
Since 
$
{\s L}D=0\Leftrightarrow D={\s S}D'$ is equivalent to
\be
D^*{\s L}^*=0\Leftrightarrow D^*=D'^*{\s S}^*
\ee
the parametrizability condition (\ref{ap})
in the case $A={\s D}'$  is
equivalent to the existence of free resolution (\ref{res2}).
\end{remark}

Suppose that the free resolution 
(\ref{res2}) exists. Let ${\s S}:={\s S}^{**}$.
As above we see that the condition $\mathrm{Ext}^1(N,{\s D}')=0$
is equivalent to
compatibility condition
\be                        \label{ab}
{\s L}^*\tilde D=\tilde D'\Leftrightarrow {\s S}^*\tilde D'=0.
\ee
Let $D=\tilde D^*,\ D'=\tilde D'^*$.
Using the relation $ ({\s L}^*\tilde D)^*
=D{\s L}=.{\s L}D$ and similarly        $({\s S}^*\tilde D')^*
=.{\s S}D'$
one sees that the compatibility condition (\ref{ab})
is equivalent to
\be                        \label{ab1}
.{\s L}D=D'\Leftrightarrow .{\s S}D'=0.
\ee
Hence the condition  
\be           \label{exa}
{\rm Ext}^1(N,{\s D}')=0
\ee
is equivalent to the assumption (\ref{eq:4.2.1}) of Theorem
\ref{th:4.2.1}.

\begin{example}
A.
Let $n=2$. Denote $\partial_1={\partial\over{\partial x_1}},\
\partial_2={\partial\over{\partial x_2}}$.
We choose
\be
{\s D}'=\left\{\qmatrix{P(\partial_1,\partial_2)&0\cr 0&0\cr}|\
P(\partial_1,\partial_2)\ {is\ PDO\ with\ constant\ coefficients}\right\}.
\ee

Consider the following
system (without boundary values)
\be
\partial_1v_2-\partial_2v_1=0.
\ee
In this case we have
\be
{\s L}=\qmatrix{{\s L}_{11}& {\s L}_{12}\cr},\ u=\qmatrix{u_1\cr u_2\cr},\
u_1=\qmatrix{v_1\cr w_1\cr},\ u_2=\qmatrix{v_2\cr w_2\cr}
\ee
where ${\s L}_{11} =\qmatrix{\partial_1&0\cr 0&0\cr},\
{\s L}_{12} =\qmatrix{-\partial_2&0\cr 0&0\cr}.$ We find
that $\ol N=2,\ \ol N_1=1$.

B. The formal adjoint exists and it is
\be
{\s L}^*=
\qmatrix{{\s L}_{11}^*\cr {\s L}_{12}^*\cr}
\ee
where
${\s L}_{11}^* =\qmatrix{-\partial_1&0\cr 0&0\cr},\
{\s L}_{12}^* =\qmatrix{\partial_2&0\cr 0&0\cr}.$ The module $N$ is given by
\be
N={\s D}'/{\s D}'^2{\s L}^*.
\ee
Consider the existence of the resolution (\ref{res2}).
Let $(D',D'')\in {\s D}'^2$. Denote
$D'=\qmatrix{P'&0\cr 0&0\cr}$, $D''=\qmatrix{P''&0\cr 0&0\cr}$
where
$P'=P'(\partial_1,\partial_2),\
P''=P''(\partial_1,\partial_2)$.
We find that
\be
.{\s L}^*(D',D'')=0\Leftrightarrow \qmatrix{D'&D''\cr}{\s L}^*=0
\ee
if and only if
\be
\qmatrix{P'&0\cr 0&0\cr}\qmatrix{-\partial_1&0\cr 0&0\cr}
+\qmatrix{P''&0\cr 0&0\cr}\qmatrix{\partial_2&0\cr 0&0\cr}=0
\Leftrightarrow
-P'\partial_1+P''\partial_2=0.\nonumber
\ee
Furthermore, the symbol of $-P'\partial_1+P''\partial_2$ is the
polynomial
$-P'(\xi_1\xi_2)\xi_1+P''(\xi_1,\xi_2)\xi_2$ and one knows from
algebra that there exists a polynomial $Q$ such that   (\cite{lang93}, p.
113)
\be
-P'(\xi_1\xi_2)\xi_1+P''(\xi_1,\xi_2)\xi_2=0\Leftrightarrow
P'(\xi_1\xi_2)=Q(\xi_1,\xi_2)\xi_2,\ P''(\xi_1,\xi_2)
=Q(\xi_1,\xi_2)\xi_1 .
\ee
Hence    the condition $-P'\partial_1+P''\partial_2=0$ is equivalent
that
\be
P'=Q(\partial_1,\partial_2)\partial_2,\
P''=Q(\partial_1,\partial_2)\partial_1 .
\ee
Denote ${\s S}^*
=
\qmatrix{\qmatrix{\partial_2&0\cr
0&0\cr}&  \qmatrix{\partial_1&0\cr 0&0\cr}\cr}.$
As a conclusion we see that
\bea
&&
.{\s L}^*(D',D'')=0\Leftrightarrow
(D',D'')= Q(\partial_1,\partial_2)\qmatrix{\qmatrix{\partial_2&0\cr
0&0\cr}&  \qmatrix{\partial_1&0\cr 0&0\cr}\cr}=Q(\partial_1,\partial_2)
{\s S}^*\nonumber
\\
&&
\Leftrightarrow (D',D'')= .{\s S}^*Q.
\eea
Hence $\ol N'=1$ and the free resolution (\ref{res2}) with finitely
generated modules exists.

C. It is easy to see that $\mathrm{Ext}^1(N,{\s D}')=0$.
Furthermore,
\be
{\s S}={\s S}^{**} =
-\qmatrix{\qmatrix{\partial_2&0\cr 0&0\cr}\cr
\qmatrix{\partial_1&0\cr 0&0\cr}\cr}.
\ee
$\mathrm{Ext}^1(\tilde N,{\s R})=0$ if and only if
\be
{\s S}U=\qmatrix{u_1\cr u_2\cr }   \Leftrightarrow {\s L}
\qmatrix{u_1\cr u_2\cr }=0
\ee
where $U=\qmatrix{V_1\cr W_1\cr}\in {\s R}$ and
$u_1=\qmatrix{v_1\cr w_1\cr},\ u_2=\qmatrix{v_2\cr w_2\cr}\in {\s R}$.
Hence we see that
$\mathrm{Ext}^1(\tilde N,{\s R})=0$ if and only if
\be
\nabla V_1=(v_1,v_2)   \Leftrightarrow \partial_1v_2- \partial_2v_1=0
\ee
which is valid if the set $G$ is convex, for example.

\end{example}

\subsection{Projectivity}\label{sec4.4}

We finally shortly treat  the projectiveness of the underlying modules.
The module $M$ is projective if and only if
$\mathrm{Ext}^1(M,A)=0$ for any left module $A$. The {\it projective
dimension}
is defined by
\be
{\rm P-dim}\ M=\inf_{n\geq 0}\{{\rm Ext}^{n+1}(M,A)=0
\ {\rm for\ any}\ A\}.
\ee
We see that   { P-dim} $M=0$ if and only if $M$ is projective.
For non-projective modules $M$ projective dimension measures
how far $M$ is from being projective.

The following theorem gives a  tool to show the
projectiveness in our case. 

\begin{theorem}\label{algp}
The module $M$ is projective  if and only if
 there exists a ${\s D}'$-homomorphism $.{\s P}:{\s D}'^{\ol N}\to
{\s D}'^{\ol N_1}$ (so called {\rm lift}) such that
\be\label{4.6a}
.{\s L}\circ .{\s P}\circ .{\s L}=.{\s L}.
\ee
\end{theorem}

The proof of  theorem can be found e.g. in  \cite{pommaret02} .
The condition (\ref{4.6a}) is equivalent to the existence of
a matrix ${\s P}\in {\s D}'(\ol N\times\ol N_1)$ such that
\be
{\s L}{\s P}{\s L}={\s L}.
\ee

\begin{corollary}\label{algpc}
Suppose that
 there exists a ${\s D}'$-homomorphism $.{\s P}:{\s D}'^{\ol N_1}\to
{\s D}'^{\ol N}$ (so called {\rm left inverse}) such that
\be\label{4.6aa}
.{\s P}\circ .{\s L} ={\s I}_{\ol N_1}
\ee
where ${\s I}_{\ol N_1}$ is the identity in ${\s D}'^{\ol N_1}$.
Then the module $M$ is projective.
\end{corollary}

The condition (\ref{4.6aa}) is equivalent to the existence of
a matrix ${\s P}\in {\s D}'(\ol N\times\ol N_1)$ such that
\be
{\s L}{\s P}={\s I}_{\ol N_1}:=\qmatrix{{\s I}& 0\cdots 0\cr
\vdots&\vdots&\vdots&\vdots\cr 0&0&\cdots&{\s I}\cr}\in
{\s D}'(\ol N_1\times\ol N_1).
\ee

Projectivity gives the following sufficient condition for the
parametrizability
\begin{theorem}\label{th10}
Suppose that there exists
a matrix ${\s P}\in {\s D}'(\ol N\times\ol N_1)$ such that
\be   \label{dec}
{\s L}{\s P}{\s L}={\s L}.
\ee
Then
the control system
${\s L}X=0$ is ${\s D}'$-parametrizable in ${ A}^{\ol N}$.
\end{theorem}

\begin{proof}
From (\ref{dec}) we get
\be
{\s L} ({\s I}_{\ol N}-{\s P}{\s L})=0.
\ee
Hence the operator ${\s S}= {\s I}_{\ol N}-{\s P} {\s L}$
satisfies the following condition: If $X={\s S}X'$ then ${\s L}X=0$.
Conversely, suppose that ${\s L}X=0$. Then
\be
X=({\s I}_{\ol N}-{\s P} {\s L})X+({\s P} {\s L})X=
{\s S}X
\ee
which completes the proof.
\end{proof}

A special case of the Theorem \ref{th10} is

 \begin{corollary}\label{co10}
Suppose that the matrix 
$ {\s L}$ 
 has a {\rm right inverse} that is, 
there exists a matrix ${\s P}\in {\s D}(\ol N\times\ol N_1)$
such that
\be             \label{deco}
 {\s L} {\s P}={\s I}_{\ol N_1}.
\ee
Then 
the control system ${\s L}X=0$ is ${\s D}'$-parametrizable
 in ${ A}^{\ol N}$.
\end{corollary}

Under the assumptions of Theorem
\ref{th10} the parametrization is given by
\be\label{pka}
{\s S}X'= ({\s I}_{\ol N}-{\s P}{\s L})X',\ X'\in A^{\ol N}.
\ee
We find that the projectiveness is a structural property
under which the parametrizability does not depend on
the module ${A}$.

\begin{remark}
A.
The module $M$ is {\rm flat} if and only if Tor-functor is exact.
Projective  module is flat but there exist nonprojective flat
modules. 
Flatness is a useful additional concept in the study
of these system modules.

B.
For example, for integral domains
${\s D}'$ one sees by a direct computation
that (\ref{ab1})
or (\ref{eq:4.2.1}) is equivalent to the  {\rm torsion freeness} of
$M$ that is,
$\tilde D[D]=0,\ \tilde D\in {\s D}',\ [D]\in M$ 
if and only if $\tilde D=0$ or $[D]=0$. 
For nonintegral domains this need not be valid.
Recall that
\be
\textrm{free}\subset \textrm{stably\ free}\subset \textrm{projective}\subset
\textrm{flat}\subset \textrm{torsion free}.
\ee

\end{remark}

\begin{example}      \label{ex:4.4.1}

Reconsider the Example \ref{ex:4.2.1}.
Recall that
\be
{\s L}:=
\qmatrix{{\s L}_1&{\s L}_2\cr}=
\qmatrix{\qmatrix{r^+(\partial_2^2-\partial_1^2)&0\cr r'I&0\cr}
&\qmatrix{-I&0\cr 0&0\cr}\cr}.
\ee
The corresponding complex related to the module
$M$ is
\be
\xymatrix{
\ar[r] & {\s D}\ar[r]^{.{\s L}} & {\s D}^2\ar[r]^{\pi} & M\ar[r] & 0.
}
\ee
Furthermore,  $D.{\s L}=0$ for ${D}=\qmatrix{{\s A}_1&{\s A}_2
\cr {\s A}_3&{\s A}_4\cr}\in
{\s D}$
if and only if
\be
{\s A}_1={\s A}_3={\s A}_2r'I={\s A}_4r'I=0.
\ee
The last two equations imply that also
${\s A}_2={\s A}_4=0$ and so $D=0$. Hence
we get the free resolution
\be                            \label{tu1}
\xymatrix{
0  \ar[r] & {\s D} \ar[r] & {\s D}^2 \ar[r]^{.{\s L}} & M \ar[r]^{\pi} &  0.
}
\ee

The truncated dual sequence is
\be
\begin{xy}
\xymatrix{
0  \ar[r] & {{\rm Hom}({\s D}^2,{\s D})}
\ar[r]^{\Phi} & {\rm Hom}({\s D},{\s D})
\ar[r] & 0
}
\end{xy}
\ee
where $\Phi(\phi)={\rm Hom}(.{\s L},{\s D})
\phi=\phi\circ .{\s L}$. The homomorphism $\Phi$ is
surjective:  
Due to Lemma \ref{le1} we see that $\Phi(\phi)=\psi$
for $\psi\in {\rm Hom}({\s D},{\s D})$ if and only if
\be                                   \label{tu}
\sum_{j=1}^2{\s L}_j\phi(E_j)=\psi(E_1).
\ee
Choose $\phi(D_1,D_2)=(D_1,D_2){\s P}\psi(E_1)$ where
\ben
{\s P}=
\qmatrix{ \qmatrix{0&K\cr 0&0\cr}\cr
\qmatrix{-I&r^+(\partial_2^2-\partial_1^2)K\cr 0&0\cr}
\cr}.
\een
Here, as above, $K$
is the potential operator
\be
Kg(x,t)=  xg(1,t)+(1-x)g(0,t)
\ee
(observe that $r'Kg=g$).
Then (\ref{tu})
holds and so $\Phi$ is surjective. Hence $\mathrm{Ext}^1(M,{\s D})=0$.

Actually  $M$ is projective.
The sequence (\ref{tu1}) splits and the 
lift can be chosen to be $.{\s P}$.
 We can compute the parametrization  by
(\ref{pka})
and the result is
\bea\label{sec3.5:eq1}
{\s S}&=&\qmatrix{\qmatrix{I&0
\cr 0&I'\cr}&\qmatrix{0&0\cr 0&0\cr}\cr
\qmatrix{0&0\cr 0&0\cr}&\qmatrix{I&0\cr 0&I'\cr}\cr}- 
{\s P}{\s L}\\
&=&
\qmatrix{ \qmatrix{I-Kr'I
&0\cr 0&I'\cr}&\qmatrix{0&0\cr 0&0\cr}\cr
\qmatrix{(r^+(\partial_2^2-\partial_1^2)
(-I+Kr'I)&0\cr 0&0\cr}&
\qmatrix{0&0\cr 0&I'\cr}\cr}. \nonumber
\eea
By (\ref{sec3.5:eq1}) the parametrization
in the sense of the subsection \ref{sec3.1}
is given elementwise (that is in the case $A={\s R}$) by 
\bea       \label{pa1}
v_1&=&f_1-Kr'f_1
\nonumber\\
v_2 &=& 
r^+(\partial_2^2-\partial_1^2)f_1-
r^+(\partial_2^2-\partial_1^2)Kr'f_1,\ f_1\in C^\infty(\ol G\times\Delta) 
\nonumber\\
w_1&=& g_1,\ g_1\in C^\infty(\partial G\times\Delta).
\eea
Note that $Kr'f_1= xf_1(1,\cdot)+(1-x)f_1(0,\cdot)$. The last equation
of (\ref{pa1}) is due to the notational covensions and it is
superflous.
\end{example}

\begin{example}\label{ex:4.4.2}
Let $n=2$ and let $\Delta v_1={\q {v_1}{x_1}}+{\q {v_1}{x_2}}$.
Consider the system
\bea\label{eq:4.4.1}
&&
(1-\Delta) v_1=0\\
&&
{\p {v_1}\nu}|_{\partial G}=w_1.\nonumber\\
\eea
Here we have $N=m=N_1=m_1=1$ and so $\ol N=\ol N_1=1$. In addition,
\be
{\s L}=\qmatrix{1-\Delta&0\cr r'{\partial\over{\partial \nu}}&-I'\cr}
\ee
where $I'$ is the identity operator on $C^\infty(\partial G)$.
Note that $r'{\partial\over{\partial \nu}}=
r'(\nu_1{\partial\over{\partial x_1}}
+\nu_2{\partial\over{\partial x_2}})$.
We seek a matrix ${\s P}\in {\s D}$ such that
\be  \label{eq:4.4.2}
{\s L}{\s P}={\s I}_1=\qmatrix{I&0\cr 0&I'\cr}.
\ee

Let ${\s P}=\qmatrix{{\s P}_1&{\s P}_2\cr {\s P}_3&{\s P}_4\cr}$.
We find that the  condition (\ref{eq:4.4.2}) is equivalent
to
\be\label{eq:4.4.3}
(1-\Delta){\s P}_1=I,\ (1-\Delta) {\s P}_2=0,\
r'{\partial\over{\partial \nu}}
{\s P}_1-{\s P}_3=0,\
r'{\partial\over{\partial \nu}}
{\s P}_2-{\s P}_4=I'.
\ee

We can choose ${\s P}_2=0$ and then ${\s P}_4=-I'$.
Let  ${\s P}_1=r^+A$ where $r^+A$ is the pseudo-differential operator
with symbol $1/(1+||\xi||^2)$. Then $(1-\Delta) {\s P}_1=I$.
From
$r'{\partial\over{\partial \nu}}
{\s P}_1-{\s P}_3=0$ we get
${\s P}_3=r'{\partial\over{\partial \nu}}
{\s P}_1$. Hence     we can choose
\be\label{eq:4.4.4}
{\s P}=\qmatrix{r^+A&0\cr  r'{\partial\over{\partial \nu}}(r^+A)&
-I'\cr}.
\ee
The corresponding parametrization is
\bea
&&
{\s S}={\s I}_1-{\s P}{\s L}
=
\qmatrix{I&0\cr 0&I'\cr}-
\qmatrix{r^+A&0\cr  r'{\partial\over{\partial \nu}}(r^+A)&
-I'\cr}
\qmatrix{1-\Delta&0\cr r'{\partial\over{\partial \nu}}&-I'\cr}
\\
&&
=
\qmatrix{I-(r^+A)(1-\Delta)&0\cr
-r'{\partial\over{\partial \nu}}(r^+A)(1-\Delta)
+r'{\partial\over{\partial \nu}}&
0\cr}.
\eea
Elementwise the parametrization is given by
\bea
&&
v_1=    (I-(r^+A)(1-\Delta))f_1\\
&&
w_1=
-r'{\partial\over{\partial \nu}}(r^+A)(1-\Delta) f_1
+r'{\p {f_1}\nu}.\nonumber
\eea

The matrix ${\s P}$ can be chosen several other ways. For example,
we can choose ${\s P}_1$ such that $(1-\Delta){\s P}_1=I,
\ r'{\partial\over{\partial \nu}}{\s P}_1=0$.
Then ${\s P}_1=r^+A+B$ where $B$ is a singular Green operator
(\cite{grubb91}).
Choosing ${\s P}_2=0$ we get
\be\label{eq:4.4.4a}
{\s P}=\qmatrix{r^+A+B&0\cr 0&
-I'\cr}.
\ee
Parametrization can be computed and it is of the form
\bea
&&
v_1=    (I-(r^+A+B)(1-\Delta))f_1\\
&&
w_1=r'{\p {f_1}\nu}.\nonumber
\eea

\end{example}

\end{document}